\magnification=\magstephalf
\newdimen\fullhsize \newdimen\hstitle \newdimen\hsbody
\tolerance=1000\hfuzz=2pt
%
\baselineskip=14pt plus 2pt minus 1pt
\parskip=4pt
\parindent=.31truein
\hoffset=.6truein
\voffset=.8truein
\hsbody=\hsize \hstitle=\hsize 
\hsize=12.6truecm
\vsize=19.8truecm
\catcode`\@=11 
\newcount\yearltd\yearltd=\year\advance\yearltd by -1900

%
%
\def\draftmode{\def\draftdate{{\rm preliminary draft:
\number\month/\number\day/\number\yearltd\ \ \hourmin}}%
\headline={\hfil\draftdate}\writelabels\baselineskip=14pt plus 2pt minus 2pt
{\count255=\time\divide\count255 by 60 \xdef\hourmin{\number\count255}
	\multiply\count255 by-60\advance\count255 by\time
   \xdef\hourmin{\hourmin:\ifnum\count255<10 0\fi\the\count255}}}

\def\draftdate{}
\def\nolabels{\def\eqnlabel##1{}\def\eqlabel##1{}\def\reflabel##1{}}
\def\writelabels{\def\eqnlabel##1{%
{\escapechar=` \hfill\rlap{\hskip.09in\string##1}}}%
\def\eqlabel##1{{\escapechar=` \rlap{\hskip.09in\string##1}}}%
\def\reflabel##1{\noexpand\llap{\string\string\string##1\hskip.31in}}}
\nolabels
%
\global\newcount\secno \global\secno=0
\global\newcount\meqno \global\meqno=1
\def\appendix#1#2{\global\meqno=1\xdef\secsym{\hbox{#1.}}\bigbreak\bigskip
\noindent{\bf Appendix #1. #2}\par\nobreak\medskip\nobreak}
%
%
\def\eqnn#1{\xdef #1{(\secsym\the\meqno)}%
\global\advance\meqno by1\eqnlabel#1}
\def\eqna#1{\xdef #1##1{\hbox{$(\secsym\the\meqno##1)$}}%
\global\advance\meqno by1\eqnlabel{#1$\{\}$}}
\def\eqn#1#2{\xdef #1{(\secsym\the\meqno)}\global\advance\meqno by1%
$$#2\eqno#1\eqlabel#1$$}
%
%
\global\newcount\refno \global\refno=1
\newwrite\rfile
\def\ref{[\the\refno]\nref}
\def\nref#1{\xdef#1{[\the\refno]}\ifnum\refno=1\immediate
\openout\rfile=refs.tmp\fi\global\advance\refno by1\chardef\wfile=\rfile
\immediate\write\rfile{\noexpand\item{#1\ }\reflabel{#1}\pctsign}\findarg}
\def\findarg#1#{\begingroup\obeylines\newlinechar=`\^^M\pass@rg}
{\obeylines\gdef\pass@rg#1{\writ@line\relax #1^^M\hbox{}^^M}%
\gdef\writ@line#1^^M{\expandafter\toks0\expandafter{\striprel@x #1}%
\edef\next{\the\toks0}\ifx\next\em@rk\let\next=\endgroup\else\ifx\next\empty%
\else\immediate\write\wfile{\the\toks0}\fi\let\next=\writ@line\fi\next\relax}}
\def\striprel@x#1{} \def\em@rk{\hbox{}} {\catcode`\%=12\xdef\pctsign{
\def\semi{;\hfil\break}
\def\addref#1{\immediate\write\rfile{\noexpand\item{}#1}} 
\def\listrefs{\vfill\eject\immediate\closeout\rfile
\baselineskip=18pt\centerline{{\bf References}}\bigskip{\frenchspacing%
\escapechar=` \input refs.tmp\vfill\eject}\nonfrenchspacing}
\def\startrefs#1{\immediate\openout\rfile=refs.tmp\refno=#1}
%
%
\font\titlerm=cmbx12 scaled\magstep5\font\titlerms=cmbx10 scaled\magstep5
\font\titlermss=cmbx7 scaled\magstep3 \font\titlei=cmmi10 scaled\magstep3
\font\titleis=cmmi7 scaled\magstep4 \font\titleiss=cmmi5 scaled\magstep4
\font\titlesy=cmsy10 scaled\magstep4 \font\titlesys=cmsy7 scaled\magstep4
\font\titlesyss=cmsy5 scaled\magstep4 \font\titleit=cmsl10 scaled\magstep4
\skewchar\titlei='177 \skewchar\titleis='177 \skewchar\titleiss='177
\skewchar\titlesy='60 \skewchar\titlesys='60 \skewchar\titlesyss='60
\font\cap=cmr17 scaled\magstep2\font\caps=cmr12 scaled\magstep2
\def\titlefont{\def\rm{\fam0\titlerm}
\textfont0=\titlerm \scriptfont0=\titlerms \scriptscriptfont0=\titlermss
\textfont1=\titlei \scriptfont1=\titleis \scriptscriptfont1=\titleiss
\textfont2=\titlesy \scriptfont2=\titlesys \scriptscriptfont2=\titlesyss
\textfont\itfam=\titleit \def\it{\fam\itfam\titleit} \rm}
\def\abstractfont{\tenpoint}

\def\tenpoint{\def\rm{\fam0\tenrm}
\textfont0=\tenrm \scriptfont0=\sevenrm \scriptscriptfont0=\fiverm
\textfont1=\teni  \scriptfont1=\seveni  \scriptscriptfont1=\fivei
\textfont2=\tensy \scriptfont2=\sevensy \scriptscriptfont2=\fivesy
\textfont\itfam=\tenit \def\it{\fam\itfam\tenit}
\textfont\bffam=\tenbf \def\bf{\fam\bffam\tenbf} \rm}
%
%
\def\newcap#1#2{
\global\ftno=0
\global\advance\secno by1
\xdef\secsym{\the\secno.}\global\meqno=1
\bigbreak
\bigskip\bigskip
\centerline{\cap CAPITULO #1}
\vskip.55in
\centerline{\bf\titlefont #2}
\vskip 3mm}
\xdef\secsym{}
%
%
\font\arm=cmbx12 scaled\magstep2 \font\arms=cmbx10 scaled\magstep2
\font\armss=cmbx7 scaled\magstep2 \font\ai=cmmi10 scaled\magstep2
\font\ais=cmmi7 scaled\magstep3 \font\aiss=cmmi5 scaled\magstep3
\font\asy=cmsy10 scaled\magstep3 \font\asys=cmsy7 scaled\magstep3
\font\asyss=cmsy5 scaled\magstep3 \font\ait=cmsl10 scaled\magstep3
\skewchar\ai='177 \skewchar\titleis='177 \skewchar\aiss='177
\skewchar\asy='60 \skewchar\titlesys='60 \skewchar\asyss='60
\def\afont{\def\rm{\fam0\arm}
\textfont0=\arm \scriptfont0=\arms \scriptscriptfont0=\armss
\textfont1=\ai \scriptfont1=\ais \scriptscriptfont1=\aiss
\textfont2=\asy \scriptfont2=\asys \scriptscriptfont2=\asyss
\textfont\itfam=\ait \def\it{\fam\itfam\ait} \rm}
%
%
\def\title#1{
\vskip 0.4cm
\noindent
\centerline{\bf #1}
\vskip 0.1cm}
%
%
\font\brm=cmbx12 \font\brms=cmbx10
\font\brmss=cmr10  \font\bi=cmmi8
\font\bis=cmmi5  \font\biss=cmmi5
\font\bsy=cmsy7  \font\bsys=cmsy7
\font\bsyss=cmsy7  \font\bit=cmsl10
\skewchar\bi='177 \skewchar\titleis='177 \skewchar\biss='177
\skewchar\bsy='60 \skewchar\titlesys='60 \skewchar\bsyss='60
\def\bfont{\def\rm{\fam0\brm}
\textfont0=\brm \scriptfont0=\brms \scriptscriptfont0=\brmss
\textfont1=\bi \scriptfont1=\bis \scriptscriptfont1=\biss
\textfont2=\bsy \scriptfont2=\bsys \scriptscriptfont2=\bsyss
\textfont\itfam=\bit \def\it{\fam\itfam\bit} \rm}
\def\subtitle#1{
\vskip 0.4cm
\noindent
{\bf\bfont #1}
\vskip 0.1cm}
%
%
\font\srm=cmr8 \font\srms=cmr5
\font\srmss=cmr5
\def\sfont{\def\rm{\fam0\srm}
\textfont0=\srm \scriptfont0=\srms \scriptscriptfont0=\srmss
\rm}
\def\feet#1{
\baselineskip=10pt
{\hskip -10mm\sfont #1}
}
%
%
\def\sifont{\def\it{\fam0\srm}
\textfont0=\srm \scriptfont0=\srms \scriptscriptfont0=\srmss
\it}
\def\feit#1{
{\sifont #1}
}
%
%
\def\parsl{\raise.15ex\hbox{/}\kern-.57em\partial}
\def\pasl{\partial \kern-1.45mm {/}}
\def\ii{\'{\i}}
\def\sc{\scriptscriptstyle}
\def\s{* d_a}
\def\lb{\lbrack}
\def\rb{\rbrack}
\def\O{\cal O}
\def\a{\alpha}
\def\r{\rho}
\def\e{\eta}
\def\s{\sigma}
\def\G{\gamma}
\def\lr{\lbrace}
\def\rr{\rbrace}
\def\ep{\epsilon}
\def\v{\vert}
\def\t{\theta}
\def\lf{\lfloor}
\def\rf{\rfloor}
\def\l{\lambda}
\def\L{\Lambda}
%
%
\font\ptlerm=cmr17 scaled\magstep4\font\ptlerms=cmr17 scaled\magstep4
%
%
$\,$
\overfullrule=0pt
%
%
\def\al{\alpha}
\def\be{\beta}
\def\ch{\chi}
\def\ga{\gamma}
\def\de{\delta}
\def\ep{\varepsilon}
\def\ze{\zeta}
\def\io{\iota}
\def\ka{\kappa}
\def\la{\lambda}
\def\na{\nabla}
\def\ro{\varrho}
\def\si{\sigma}
\def\om{\omega}
\def\ph{\varphi}
\def\ta{\tau}
\def\th{\theta}
\def\te{\vartheta}
\def\up{\upsilon}
\def\Ga{\Gamma}
\def\De{\Delta}
\def\La{\Lambda}
\def\Si{\Sigma}
\def\Om{\Omega}
\def\Te{\Theta}
\def\Th{\Theta}
\def\Up{\Upsilon}
\def\CO{{\cal{O}}}
\def\sign{{\rm sign}}
\def\inbar{\,\vrule height1.5ex width.4pt depth0pt}
\def\IC{\relax{\raise .7mm\hbox{${\scriptstyle\vert}$}}\hskip -1.5mm {\rm C}}
\def\IN{\relax{\raise .4mm\hbox{${\scriptscriptstyle\vert}$}}\hskip -1.4mm
{\rm N}}

\def\R{\scriptscriptstyle R}
\def\scc{\scriptscriptstyle}
\def\ba{${\scriptstyle /}$}
\def\ve{${\scriptstyle \vert}$}
\def\bas{${\scriptscriptstyle /}$}
\def\Cs{C \kern-2.2mm \raise.4ex\hbox {\ba} \kern.2mm}
\def\Css{C \kern-1.9mm \raise.3ex\hbox {\bas} \kern.8mm}
\def\Zs{Z \kern-2.2mm \raise.4ex\hbox {\ba} \kern.2mm}
\def\Rs{{\rm R} \kern-2.9mm \raise.35ex\hbox {\ve} \kern.4mm\hskip 1.3mm}
\def\Ns{{\rm N} \kern-2.9mm \raise.35ex\hbox {\ve} \kern.4mm\hskip 1.3mm}
\def\boxe{\sqcup \kern-2.3mm \sqcap}
%
\magnification=\magstephalf


\hsize = 14.6truecm
\vsize = 21.0truecm

\nopagenumbers

\headline{
           \rm \hfill\the\pageno}

\input epsf


\hoffset=6truemm
\baselineskip=13pt
\parskip 2pt plus 1pt
\parindent=15pt

\font\notefont=cmr9

\vbadness=10000
\widowpenalty=10000
\clubpenalty=10000

\font\afont=cmss12
\font\bfont=cmss9
\font\cfont=cmssi9
\font\hfont=cmcsc10
\font\pfont=cmss10
\font\tfont=cmss12 at 18pt
\font\reffont=cmcsc10
\font\secfont=cmbx12
\def\subsecfont{\bf}
\def\subsubsecfont{\sl}
\font \abscwifont=cmss10

\font\footfont=cmr8
\def\ref#1{{\reffont#1}}

\def\sq{\quad}
\def\sqq{\qquad}
\def\cl{\centerline}

\newcount\secnum
\newcount\subsecnum
\newcount\subsubsecnum
\newcount\remarknum
\newcount\lemmanum
\newcount\thmnum
\newcount\vbnum

\def\sect#1{\advance\secnum by 1\subsecnum=0 \equationnumber=0
                     \remarknum=0\lemmanum=0\thmnum=0
                  \vskip 18pt
                \leftline{\secfont\the\secnum.\    #1}
            \vskip -12pt
                \message{#1}\nobreak\noindent}

	\def\subsect#1{\advance\subsecnum by 1\subsubsecnum=0
 			           \vskip 14pt
                \leftline{\subsecfont \the\secnum.\the\subsecnum.\ #1}
          			\vskip -8pt
                \message{#1}\nobreak \noindent}

\def\subsubsect#1{\advance\subsubsecnum by 1
							\vskip 12pt

\leftline{\subsubsecfont\the\secnum.\the\subsecnum.\the\subsubsecnum. \ #1}
          			\vskip -6pt
               \message{#1}\nobreak \noindent}

\def\remark{\advance \remarknum by 1
         \vskip\baselineskip
         \noindent {\bf Remark \the\secnum.\the\remarknum.}\quad\rm
         }
\def\lemma{\advance \lemmanum by 1
         \vskip\baselineskip
         \noindent {\bf Lemma \the\secnum.\the\lemmanum.}\quad\it
         }
\def\theorem{\advance \thmnum by 1
         \vskip\baselineskip
         \noindent {\bf Theorem \the\secnum.\the\thmnum.}\quad\it
         }
\def\proof{
         \vskip\baselineskip
         \noindent {\bfProof }\quad\rm
         }
\def\vb{\advance \vbnum by 1
         \vskip\baselineskip
         \noindent {\bf Example 
\the\vbnum.}\quad
         }

\def\brg#1#2{{{\lower.4ex
\hbox{$\scriptstyle#1$}}\over
{\raise.4ex
\hbox{$\scriptstyle#2$}}}}

\def\remp
#1. #2\par{\medbreak\noindent{\bf#1.\enspace}{\rm#2}\par\medbreak}

\def\thep
#1. #2\par{\medbreak\noindent{\bf#1.\enspace}{\sl#2}\par\medbreak}

\def\prop
#1. #2\par{\medbreak\noindent{\bf#1.\enspace}{\rm#2}\par
 \rightline{\vrule height4pt width5.5pt depth2pt}\medbreak}
\def\proof{\bf \medbreak \noindent Proof. \rm}
\def\eoproof{{\unskip\nobreak\hfil\penalty50
	\hskip2em\hbox{}\nobreak\hfil\vrule height4pt width5.5pt depth2pt
	\parfillskip=0pt\finalhyphendemerits=0\medbreak}}

\def\w#1{{\sqrt#1}\,}
\def\kd{\partial}
\def\sq{\quad}
\def\sqs{\qquad}
\def\n{\eqno}
\def\cl{\centerline}
\def\el{\eqalign}

\def\iy{\infty}

\def\bo{{\cal O}}

\def\pd#1#2{{{\kd#1}\over{\kd#2}}}
\def\pdt#1#2{{{\kd^2#1}\over{\kd #2^2}}}

\def\br#1#2{{{#1}\over{#2}}}

\def\intp{\int_0^\iy}
\def\intr{\int_{-\iy}^\iy}

\def\RR{{{\rm I}\!{\rm R}}}
\def\NN{{{\rm I}\!{\rm N}}}
\def\RRP{{{\rm I}\!{\rm R^+}}}
\def\RRN{{{\rm I}\!{\rm R_0}}}
\def\NNP{{{\rm I}\!{\rm N^+}}}
\def\ZZ{{\hbox{Z}\!\!\hbox{Z}}}
\def\KK{{{\rm I}\!{\rm K}}}
\def\one{{{\rm 1}\hskip-0.55ex{\rm I}}}
\def\CC{\hbox{\rlap{$\,\,
  $\hbox{\vrule height6.2pt width.35pt depth-0.1pt}}$\rm C$}}
\def\QQ{\hbox{\rlap{$\,\,
  $\hbox{\vrule height6pt width.35pt depth0.1pt}}$\rm Q$}}
\def\P{\cal P}

\def\frac#1/#2{\leavevmode\kern.1em\raise.5ex\hbox{\the\scriptfont0
#1}\kern-.1em/\kern-.15em\lower.25ex\hbox{\the\scriptfont0
#2}}

%
\newcount\equationnumber	\equationnumber=0
\def\eqnum{\relax
	\global\advance\equationnumber by 1
	\equationnumberformat{\the\equationnumber}%
	}%
\def\eqname#1{\relax
	\count255=\equationnumber
	\assignnumber{EN#1}\equationnumber
	\global\equationnumber=\count255
	\global\advance\equationnumber by 1
	\ifnum\csname EN#1\endcsname=\equationnumber
	\else
		\message{The equation number for ``#1'' is incorrect!}%
	\fi
	\equationnumberformat{\csname EN#1\endcsname}%
	}%
\def\equationnumberformat#1{\eqno(\the\secnum.\equationnumbertype{#1})}%
\def\equationnumbertype#1{\number#1\relax}%
\def\referenceequation#1{\relax
	\assignnumber{EN#1}\equationnumber
	\equationnumbertype{\csname EN#1\endcsname}%
	}%
\def\forwardreferenceequation#1#2{\relax
	\global\advance\equationnumber by #2
	\assignnumber{EN#1}\equationnumber
	\global\advance\equationnumber by -1
	\global\advance\equationnumber by -#2
	\referenceequation{#1}%
	}%
%
\def\assignnumber#1#2{\relax
	\ifnum0<0\csname#1\endcsname
	\else
		\global\advance#2 by 1
		\expandafter\expandafter\expandafter
			\xdef\csname#1\endcsname{\the#2}%
	\fi
	}%

\def\fre{\forwardreferenceequation}
\def\en{\eqname}
\def\req#1{(\the\secnum.{\referenceequation{#1}})}

\def\arcsinh{{\rm arcsinh}}
\def\arccosh{{\rm arccosh}}
\def\arctanh{{\rm arctanh}}
\def\arccoth{{\rm arccoth}}

\def\erf{{\rm erf}}
\def\erfc{{\rm erfc}}

\def\phase{{\rm phase}}
\def\sign{{\rm sign}}

\def\wt{{\sqrt{2}}}

\def\phih{{\widehat \phi}}
\def\psih{{\widehat \psi}}

\def\Ai{{{\rm Ai}}}
\def\Bi{{{\rm Bi}}}

\def\phizeta{\left(\br{4\z}{1-z^2}\right)^{1/4}}

\catcode`@=11 

\font\ninerm=cmr10 at 9pt
\font\eightrm=cmr7 at 8pt
\font\sixrm=cmr7 at 6pt

\font\ninei=cmmi10 at 9pt
\font\eighti=cmmi10 at 8pt
\font\sixi=cmmi7 at 6pt
\skewchar\ninei='177
\skewchar\eighti='177
\skewchar\sixi='177

\font\ninesy=cmsy10 at 9pt
\font\eightsy=cmsy10 at 8pt
\font\sixsy=cmsy7 at 6pt
\skewchar\ninesy='60
\skewchar\eightsy='60
\skewchar\sixsy='60

\font\eightss=cmssq8

\font\eightssi=cmssqi8

\font\ninebf=cmbx10 at 9pt
\font\eightbf=cmbx10 at 8pt
\font\sixbf=cmbx7 at 6pt

\font\ninett=cmtt10 at 9pt
\font\eighttt=cmtt10 at 8pt

\hyphenchar\tentt=-1 
\hyphenchar\ninett=-1
\hyphenchar\eighttt=-1

\font\ninesl=cmsl10 at 9pt
\font\eightsl=cmsl10 at 8pt

\font\nineit=cmti10 at 9pt
\font\eightit=cmti10 at 8pt

\font\tenu=cmu10 
\newskip\ttglue

\def\eightpoint{\def\rm{\fam0\eightrm}%
  \textfont0=\eightrm \scriptfont0=\sixrm \scriptscriptfont0=\fiverm
  \textfont1=\eighti \scriptfont1=\sixi \scriptscriptfont1=\fivei
  \textfont2=\eightsy \scriptfont2=\sixsy \scriptscriptfont2=\fivesy
  \textfont3=\tenex \scriptfont3=\tenex \scriptscriptfont3=\tenex
  \def\it{\fam\itfam\eightit}%
  \textfont\itfam=\eightit
  \def\sl{\fam\slfam\eightsl}%
  \textfont\slfam=\eightsl
  \def\bf{\fam\bffam\eightbf}%
  \textfont\bffam=\eightbf \scriptfont\bffam=\sixbf
   \scriptscriptfont\bffam=\fivebf
  \def\tt{\fam\ttfam\eighttt}%
  \textfont\ttfam=\eighttt
  \tt \ttglue=.5em plus.25em minus.15em
  \normalbaselineskip=9pt
  \def\MF{{\manual opqr}\-{\manual stuq}}%
  \let\sc=\sixrm
  \let\big=\eightbig
  \setbox\strutbox=\hbox{\vrule height7pt depth2pt width\z@}%
  \normalbaselines\rm}

\def\a{\alpha}\def\b{\beta}\def\c{\gamma}\def\d{\delta} \def\vth{\vartheta}
\def\eps{\varepsilon}\def\f{\phi}\def\k{\kappa}\def\l{\lambda}\def\m{\mu}
\def\p{\pi}\def\r{\rho}\def\s{\sigma}\def\t{\tau}\def\th{\theta}
\def\x{\xi}\def\y{\eta}\def\z{\zeta}\def\om{\omega}
\def\La{\Lambda}\def\oom{\Omega}\def\G{\Gamma}\def\D{\Delta}

\def\sn{\sum_{n=0}^\iy}
\def\sk{\sum_{k=0}^\iy}
\def \som#1#2{\sum_{{#1}}^{{#2}}}

\def\nfe{\eject}
\def\vfe{\vfil\eject}
\def\oom{\Omega}
\def\wh{\widehat}

\def\gp{$C_n^{\c}(x)$}
\def\lp{$L_n^{\a}(x)$}
\def\jp{$P_n^{(\a,\b)}(x)$}

\def\gpf{C_n^{\c}(x)}
\def\lpf{L_n^{\a}(x)}
\def\jpf{P_n^{(\a,\b)}(x)}

\def\bri{\br1{2\pi i}}
\def\C{{\cal C}}
\def\H{{\cal H}}
\def\L{{\cal L}}

\def\({\left(}
\def\){\right)}

\def\[{\left[}
\def\]{\right]}

\def\aa{^{(\a)}(x)}

\noindent
%


\nref\Chester{{\reffont C. Chester, B. Friedman, and F. Ursell},
An extension of the method of steepest descent,
Proc. Cambridge Philos. Soc., {\bf 53} (1957) 599--611.}

\nref\nicoi{{\reffont Jos\'e L. Lopez and Nico M. Temme},
Asymptotic expansions of Charlier, Laguerre and Jacobi polynomials.
Submited, 2002.}

\nref\nicoii{{\reffont Jos\'e L. Lopez and Nico M. Temme},
Two-point Taylor expansions of analytic functions.
Submited, 2002.}

\nref\Raimundas{{\reffont Raimundas Vidunas and Nico M. Temme},
Symbolic evaluation of coefficients in Airy-type asymptotic expansions.
CWI Report MAS-R0118.
Accepted for publication in
{\sl Journal of Mathematical Analysis and Applications} (2001).}

\nref\walsh{{\reffont J. L. Walsh},
Interpolation and Approximation by rational functions in the complex domain,
{\sl American Mathematical Society, Providence}, 1969.}

\nref\Whittaker{{\reffont E.T. Whittaker and G.N. Watson},
A course of modern  analysis,
{\sl Cambridge University Press, London and New York}, 1927.}

\nref\wong{{\sl R. Wong}, Asymptotic Approximations of
Integrals, {\sl Academic Press, New York}, 1989.}


\tfont
\centerline      {Two-point Taylor Expansions of Analytic Functions}
\bigskip\bigskip 
\afont
\cl{Jos\'e L. L\'opez$^1$  and Nico M. Temme$^2$}
\medskip
\cfont
\cl{$^1$
Departamento de Mat\'ematica e Inform\'atica,}
\cl{Universidad P\'ublica de Navarra,
31006-Pamplona,
Spain
}
\cl{$^2$
CWI,
P.O. Box 94079,
1090 GB Amsterdam,
The Netherlands}
\cfont
\cl{e-mail: \tt jl.lopez@unavarra.es, nicot@cwi.nl}
\medskip

\bfont
\parindent=25pt
{\pfont ABSTRACT}\par\noindent
{\narrower\noindent
Taylor expansions of analytic functions are considered with respect to two
points. Cauchy-type formulas are given for coefficients and remainders in
the expansions, and the regions of convergence are indicated. It is
explained how these expansions can be used in deriving  uniform asymptotic
expansions of integrals. The method is also used for obtaining  Laurent
expansions in two points.

\vskip 0.3cm \noindent
\cfont
2000 Mathematics Subject Classification:
\bfont
30B10, 30E20, 40A30.
\par\noindent
\cfont
Keywords \& Phrases:
\bfont
two-point Taylor expansions, Cauchy's theorem,
analytic functions,
two-point Laurent expansions,
uniform asymptotic expansions of integrals.

\par\noindent   
{\parindent 35pt
\cfont  
\item{Note:\quad}
\bfont  Work carried out under project MAS1.2 
Analysis, Asymptotics and Computing. 
This report has been accepted for publication in 
{\sl Studies in Applied Mathematics}.
\par
}  
}
\rm
\parindent=15pt
\sect{Introduction}

\noindent
In deriving uniform asymptotic expansions of a certain class of integrals
one encounters the problem of expanding a function, that is analytic in some
domain $\Omega$ of the complex plane, in two points. The first mention
of the use of such expansions in asymptotics is given in \ref\Chester,
where Airy-type expansions are derived for integrals having two
nearby (or coalescing) saddle points. This reference does not give further
details about two-point Taylor expansions, because the coefficients in the
Airy-type asymptotic expansion are derived in a different way.

To demonstrate the application in asymptotics we consider the integral
\eqn\CFU{
F_b(\omega)={1\over 2\pi i}\int _{\C} e^{\omega({1\over3} z^3-b^2z)} f(z)\,dz,
}
where $\omega$ is a large positive parameter and $b$ is a parameter that may
assume small
values. The contour starts at $\infty e^{-i\pi/3}$ and terminates at
$\infty e^{i\pi/3}$, and lies in a domain where the function $f$ is
analytic.
In particular, $f$ is analytic in a domain that contains the saddle
points $\pm b$ of the exponent in the integrand. One method for obtaining an
asymptotic expansion of $F_b(\omega)$ that holds uniformly for small values of
$b$ is based on expanding $f$ at the two saddle points:
\eqn\ftexp{
f(z)=\sum_{n=0}^\infty A_n(z^2-b^2)^n + z \sum_{n=0}^\infty
B_n(z^2-b^2)^n,
}
and substitute this expansion into \CFU.  When interchanging summation
and integration, the result is a formal expansion in two series in
terms of functions related with Airy functions.  A Maple algorithm for
obtaining the coefficients $A_n$ and $B_n$, with applications to
Airy-type expansions of parabolic cylinder functions, is given in
\ref\Raimundas.

In a future paper we use expansions like \ftexp\ in order to derive
convergent expansions for orthogonal polynomials and hypergeometric
functions that also have an asymptotic nature.  The purpose of the
present paper is to give details on the two-point Taylor expansion
\ftexp, in particular on the region of convergence and on
representations in terms of Cauchy-type integrals of coefficients and
remainders of these expansions.  Some information on this type of
expansions is also given in \ref\Whittaker, p.  149,
Exercise 24.

Without referring to applications in asymptotic analysis we include
analogous properties of two-point Laurent expansions and of another
related type, the two-point Taylor-Laurent expansion.

\sect{Two-point Taylor expansions}

\noindent
We consider the expansion \ftexp\ in a more symmetric form and give
information on
the coefficients and the remainder in the expansion.

\noindent
{\bf Theorem 1.} {\it Let $f(z)$ be an analytic function on an open set
$\Omega\subset\Cs$
and $z_1$, $z_2\in\Omega$ with $z_1\ne z_2$. Then, $f(z)$ admits the two-point
Taylor expansion
\eqn\expani{
f(z)=\sum_{n=0}^{N-1}\left\lbrack
a_n(z_1,z_2)(z-z_1)+a_n(z_2,z_1)(z-z_2)\right\rbrack
(z-z_1)^n(z-z_2)^n+r_N(z_1,z_2;z),
}
where the coefficients $a_n(z_1,z_2)$ and $a_n(z_2,z_1)$ of the expansion
are given by the
Cauchy integral
\eqn\coefi{
a_n(z_1,z_2)\equiv{1\over 2\pi i(z_2-z_1)}\int_{\cal C}{f(w)\,dw
\over (w-z_1)^n(w-z_2)^{n+1}}.
}
The remainder term $r_N(z_1,z_2;z)$ is given by the Cauchy integral
\eqn\remi{
r_N(z_1,z_2;z)\equiv{1\over 2\pi i}\int_{\cal C}{f(w)\,dw\over
(w-z_1)^N(w-z_2)^N
(w-z)}(z-z_1)^N(z-z_2)^N.
}
The contour of integration ${\cal C}$
is a simple closed loop which encircles the points $z_1$ and $z_2$ (for
$a_n$) and $z$, $z_1$ and $z_2$ (for $r_N$) in
the counterclockwise direction and is contained in $\Omega$ (see
Figure 1 (a)).

The expansion \expani\ is convergent for $z$ inside the Cassini oval (see
Figure 2)
$$
O_{z_1,z_2}\equiv\lbrace z\in\Omega, \hskip 2mm
\vert(z-z_1)(z-z_2)\vert<r\rbrace
$$
where
$$
r\equiv {\rm Inf}_{w\in\Css\setminus\Omega}\left\lbrace
\vert(w-z_1)(w-z_2)\vert\right\rbrace.
$$
 In particular, if $f(z)$ is an
entire function
$(\Omega=\Cs$), then the
expansion \expani\ converges $\forall$ $z\in\Cs$.
}

\noindent
{\bf Proof.} By Cauchy's theorem,
\eqn\cauchy{
f(z)={1\over 2\pi i}\int_{\cal C}{f(w)\,dw\over w-z},
}
where ${\cal C}$ is the contour defined above (Figure 1 (a)). Now we write
\eqn\lequal{
{1\over w-z}={z+w-z_1-z_2\over (w-z_1)(w-z_2)}{1\over 1-u},
}
where
\eqn\ui{
u\equiv {(z-z_1)(z-z_2)\over (w-z_1)(w-z_2)}.
}
Now we introduce the expansion
\eqn\expanu{
{1\over 1-u}=\sum_{n=0}^{N-1}u^n+{u^N\over 1-u}
}
in \lequal\ and this in \cauchy. After straightforward calculations we
obtain \expani-\remi.

For any $z\in O_{z_1,z_2}$, we can take a contour $\cal C$ in $\Omega$ such
that
$\vert(z-z_1)(z-z_2)\vert<\vert(w-z_1)(w-z_2)\vert$ $\forall$ $w\in{\cal
C}$ (see Figure 1 (b)).
In this contour $\vert f(w)\vert$ is bounded by some constant $C$: $\vert
f(w)\vert\le C$.
Introducing these two bounds in \remi\ we see that
$\lim_{N\to\infty}r_N(z_1,z_2;z)=0$ and
the proof follows.
\hfill $\boxe$

\bigskip
\centerline{\epsfxsize=8cm \epsfbox{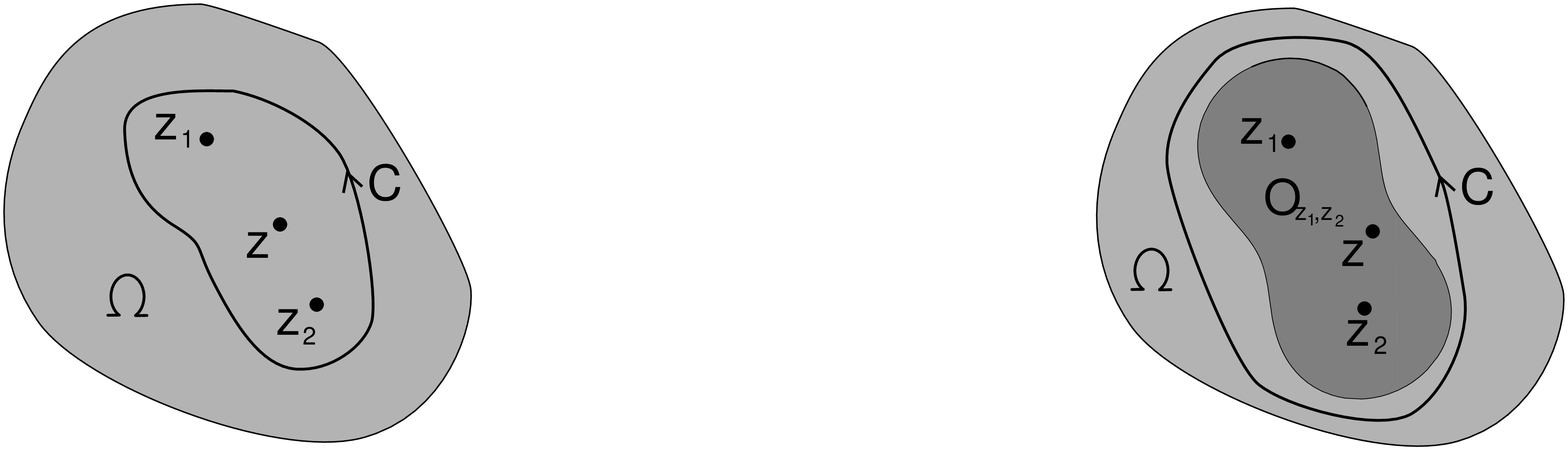}}
\centerline{(a) \hskip 6cm (b)}
\parindent=10pt

\noindent
{\bf Figure 1}. {\cfont (a) Contour $\cal C$ in the integrals
\expani-\remi. (b) For
$z\in O_{z_1,z_2}$, we can take a contour $\cal C$ in $\Omega$ which contains
$O_{z_1,z_2}$ inside and therefore,
$\vert(z-z_1)(z-z_2)\vert<\vert(w-z_1)(w-z_2)\vert$ $\forall$ $w\in{\cal C}$.}
\vskip 2mm

\rm
\parindent=15pt

\bigskip
\centerline{\epsfxsize=10cm \epsfbox{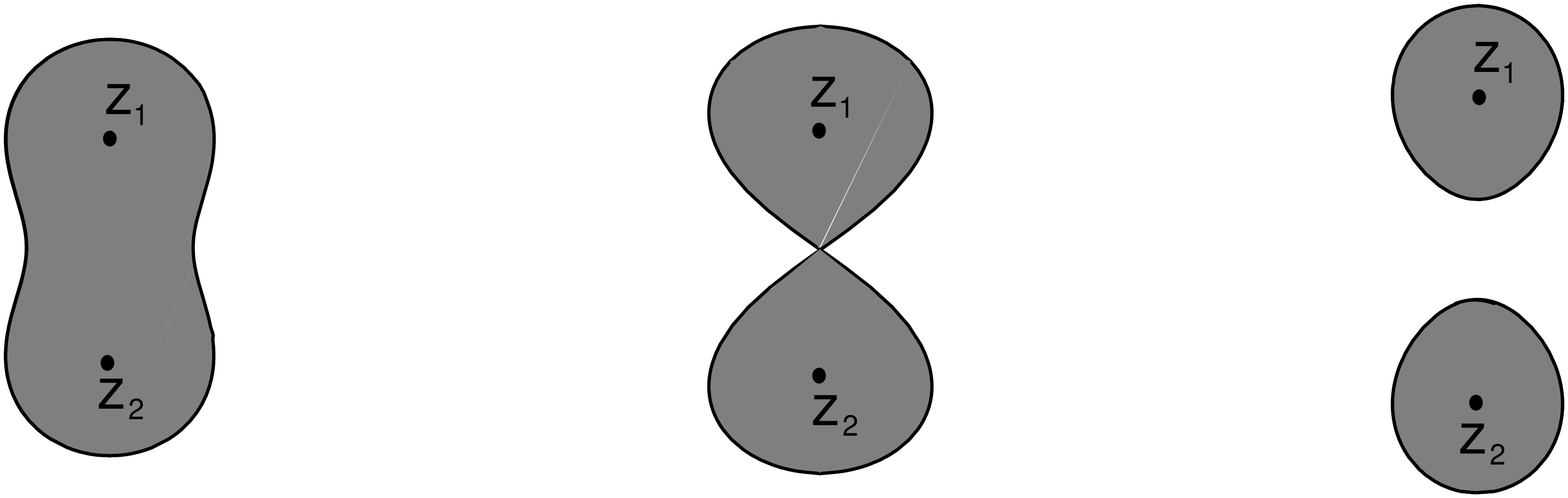}}
\centerline{(a) $4r>\vert z_1-z_2\vert^2$ \hskip 1.5cm (b) $4r=\vert
z_1-z_2\vert^2$
\hskip 1.5cm (c) $4r<\vert z_1-z_2\vert^2$}
\parindent=10pt

\noindent
{\bf Figure 2}. {\cfont Shape of the Cassini oval depending on the relative
size of the parameter $r$ and the focal distance $\vert z_1-z_2\vert$.}
\vskip 2mm

\subsect{An alternative form of the expansion}%

\smallskip
\noindent
The present expansion of $f(z)$ in the form \expani\ stresses the
symmetry of the expansion with respect to $z_1$ and $z_2$. In this
representation it is not possible, however, to let $z_1$ and $z_2$
coincide, which causes a little inconvenience (the coefficients
$a_n(z_1,z_2)$ become infinitely large as $z_1\to z_2$;
the remainder $r_N(z_1,z_2;z)$ remains well-defined).
An alternative way is the representation (cf. \ftexp),
$$
f(z)=\sum_{n=0}^{\infty}\left\lbrack
A_n(z_1,z_2)+B_n(z_1,z_2)\,z\right\rbrack
(z-z_1)^n(z-z_2)^n,
$$
and we have the  relations
$$
\eqalign{
A_n(z_1,z_2)&=-z_1a_n(z_1,z_2)-z_2a_n(z_2,z_1),\cr
B_n(z_1,z_2)&= a_n(z_1,z_2)+a_n(z_2,z_1),\cr
}
$$
which are regular when $z_1\to z_2$. In fact we have
$$
\eqalign{
A_n(z_1,z_2)&={1\over 2\pi i}\int_{\cal C}{w-z_1-z_2\over
[(w-z_1)(w-z_2)]^{n+1}} f(w)\,dw,\cr
B_n(z_1,z_2)&={1\over 2\pi i}\int_{\cal C}{f(w)\,dw\over
[(w-z_1)(w-z_2)]^{n+1}}.\cr
}
$$
Letting $z_1\to0$ and $z_2\to0$, we obtain the standard Maclaurin
series of $f(z)$ with even part (the $A_n$ series) and odd part
(the $B_n$ series).

\rm
\parindent=15pt

\subsect{Explicit forms of the coefficients}%

\smallskip
\noindent
Definition \coefi\ is not appropriate for numerical computations.
A more practical formula to compute the coefficients of the above
two-point Taylor expansion is given in the following proposition.

\noindent
{\bf Proposition 1.} {\it Coefficients $a_n(z_1,z_2)$ in the expansion
\expani\ are
also given by the formulas:
\eqn\ao{
a_0(z_1,z_2)={f(z_2)\over z_2-z_1}
}
and,  for $n=1,2,3,...$,
\eqn\coefibis{
a_n(z_1,z_2)=\sum_{k=0}^n{(n+k-1)!\over
k!(n-k)!}{(-1)^{n+1}nf^{(n-k)}(z_2)+(-1)^kkf^{(n-k)}(z_1)\over
n!(z_1-z_2)^{n+k+1}}.
}
%

\noindent
{\bf Proof.} We deform the contour of integration  ${\cal C}$ in equation
\coefi\ to any
contour  of the form ${\cal C}_1\cup{\cal C}_2$ also contained in $\Omega$,
where
${\cal C}_1$ (${\cal C}_2$) is a simple closed loop which encircles the
point $z_1$ ($z_2$) in the
counterclockwise direction and does not contain the point $z_2$ ($z_1$)
inside (see Figure 3 (a)). Then,
$$
\eqalign{
a_n(z_1,z_2)= & {1\over 2\pi i(z_2-z_1)}\left\lbrace\int_{{\cal C}_1}{f(w)\over
(w-z_2)^{n+1}}{dw\over (w-z_1)^n}+\int_{{\cal C}_2}{f(w)\over
(w-z_1)^{n}}{dw\over (w-z_2)^{n+1}}\right\rbrace= \cr &
{1\over (z_2-z_1)}\left\lbrace{1\over (n-1)!}{d^{n-1}\over dw^{n-1}}
\left.{f(w)\over(w-z_2)^{n+1}}\right\vert_{w=z_1}+
{1\over  n!}{d^{n}\over dw^{n}}
\left.{f(w)\over(w-z_1)^{n}}\right\vert_{w=z_2}\right\rbrace. \cr}
$$
From here, equations \ao-\coefibis\ follows after straightforward computations.
\hfill$\boxe$

\bigskip
\centerline{\epsfxsize=12cm \epsfbox{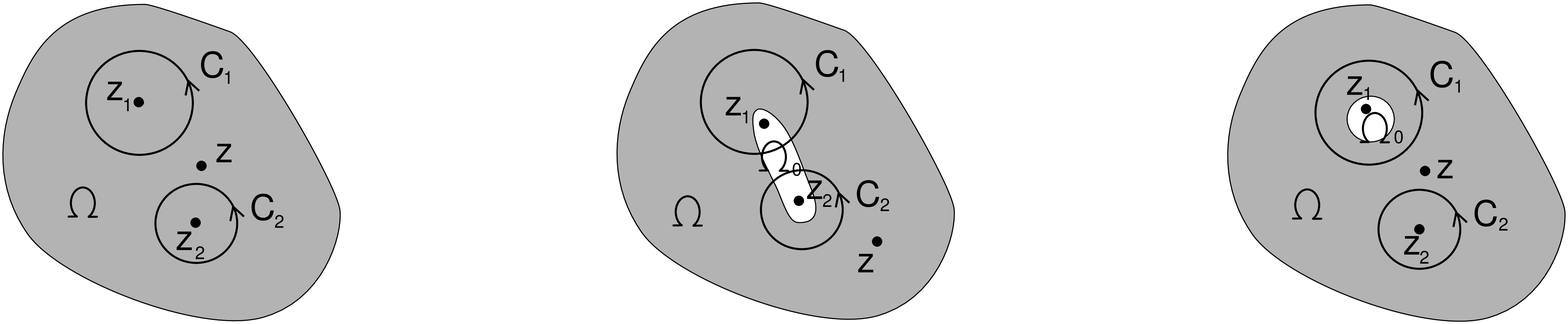}}
\centerline{(a) \hskip 4cm (b)\hskip 4cm (c)}
\parindent=10pt

\noindent
{\bf Figure 3}. {\cfont (a) The function $(w-z_2)^{-n-1}f(w)$ is analytic
inside
${\cal C}_1$, whereas $(w-z_1)^{-n}f(w)$ is analytic inside ${\cal C}_2$.
(b) The function $(w-z_2)^{-n-1}g_1(w)$ is analytic inside
${\cal C}_1$, whereas $(w-z_1)^{-n}g_2(w)$ is analytic inside ${\cal C}_2$.
(c) The function $(w-z_2)^{-n-1}g(w)$ is analytic inside
${\cal C}_1$, whereas $(w-z_1)^{-n}f(w)$ is analytic inside ${\cal C}_2$.}
\vskip 2mm
\rm
\parindent=15pt

\subsect{Two-point Taylor polynomials}%

\smallskip
\noindent
Next we can define the two-point Taylor polynomial of the function $f(z)$ at
in the following way:

\noindent
{\bf Definition 1.} {\it Let $z$ be a real or complex variable and $z_1$
and $z_2$ ($z_1\ne z_2$) two real or complex numbers.  If $f(z)$ is
$n-1-$times differentiable at those two points, we define the two-point
Taylor polynomial of $f(z)$ at $z_1$ and $z_2$ and degree $2n-1$ as
$$
P_n(z_1,z_2;z)\equiv
\sum_{k=0}^{n-1}\left\lbrack
a_k(z_1,z_2)(z-z_1)+a_k(z_2,z_1)(z-z_2)\right\rbrack
(z-z_1)^k(z-z_2)^k,
$$
where the coefficients $a_k(z_1,z_2)$ are given in \ao-\coefibis.}

\noindent
{\bf Proposition 2.} {\it In the conditions of the above definition, define
the remainder of the
approximation of $f(z)$ by $P_n(z_1,z_2;z)$ at $z_1$ and $z_2$ as
$$
r_n(z_1,z_2;z)\equiv f(z)-P_n(z_1,z_2;z).
$$
Then, (i) $r_n(z_1,z_2;z)=o(z-z_1)^{n-1}$ as $z\to z_1$ and
$r_n(z_1,z_2;z)=o(z-z_2)^{n-1}$ as $z\to z_2$.
(ii) If $f(z)$ is $n-$times differentiable at $z_1$ and/or $z_2$, then
$r_n(z_1,z_2;z)={\cal O}(z-z_1)^n$ as $z\to z_1$ and/or
$r_n(z_1,z_2;z)={\cal O}(z-z_2)^n$ as $z\to z_2$.}

\noindent
{\bf Proof.} The proof is trivial if $f(z)$ is analytic at $z_1$ and $z_2$
by using \remi.
In any case, for real or complex variable, the proof follows after
straightforward
computations by using l'H\^opital's rule and \ao-\coefibis.
\hfill $\boxe$

\noindent
{\bf  Remark 1.} {
Observe that the Taylor polynomial of $f(z)$
at $z_1$ and $z_2$
and degree $2n-1$ is the same as the Hermite's interpolation polynomial of
$f(z)$
at $z_1$ and $z_2$ with data $f(z_i)$, $f'(z_i)$,...,$f^{(n-1)}(z_i)$,
$i=1,2$.}

\sect{Two-point Laurent expansions}

\noindent
In the standard theory for Taylor and Laurent expansions much analogy
exists between the two types of expansions. For two-point expansions, we
have a similar agreement in the representations of coefficients and remainders.

\noindent
{\bf Theorem 2.} {\it Let $\Omega_0$ and $\Omega$ be closed and open sets,
respectively, of the complex plane, and $\Omega_0\subset\Omega\subset\Cs$. Let
$f(z)$ be an analytic function on $\Omega\setminus\Omega_0$
and $z_1$, $z_2\in\Omega_0$ with $z_1\ne z_2$.
Then, for any $z\in\Omega\setminus\Omega_0$,
$f(z)$ admits the two-point Laurent expansion
\eqn\expanii{
\eqalign{
f(z)= & \sum_{n=0}^{N-1}\left\lbrack
b_n(z_1,z_2)(z-z_1)+b_n(z_2,z_1)(z-z_2)\right\rbrack
(z-z_1)^n(z-z_2)^n+ \cr &
\sum_{n=0}^{N-1}\left\lbrack
c_n(z_1,z_2)(z-z_1)+c_n(z_2,z_1)(z-z_2)\right\rbrack
(z-z_1)^{-n-1}(z-z_2)^{-n-1}+ \cr & r_N(z_1,z_2;z), \cr}
}
where the coefficients $b_n(z_1,z_2)$, $b_n(z_2,z_1)$, $c_n(z_1,z_2)$ and
$c_n(z_2,z_1)$ of the expansion are given, respectively, by the  Cauchy
integrals
\eqn\coefii{
b_n(z_1,z_2)\equiv{1\over 2\pi i(z_2-z_1)}\int_{{\Gamma}_1}
{f(w)\,dw\over (w-z_1)^n(w-z_2)^{n+1}}
}
and
\eqn\coefiii{
c_n(z_1,z_2)\equiv{1\over 2\pi i(z_2-z_1)}\int_{{\Gamma}_2}
(w-z_1)^{n+1}(w-z_2)^{n}f(w)\,dw.
}
The remainder term $r_N(z_1,z_2;z)$ is given by the Cauchy integrals
\eqn\remii{
\eqalign{
r_N(z_1,z_2;z)\equiv & {1\over 2\pi i}\int_{{\Gamma}_1}{f(w)dw\over
(w-z_1)^N(w-z_2)^N
(w-z)}(z-z_1)^N(z-z_2)^N- \cr &
{1\over 2\pi i}\int_{{\Gamma}_2}{(w-z_1)^N(w-z_2)^Nf(w)dw\over
w-z}{1\over(z-z_1)^N(z-z_2)^N}.\cr}
}
In these integrals, the contours of integration ${\Gamma}_1$ and
${\Gamma}_2$ are
simple closed loops
contained in $\Omega\setminus\Omega_0$ which encircle the
points $z_1$ and $z_2$ in the counterclockwise direction.
Moreover, ${\Gamma}_2$ does not contain
the point $z$ inside,  whereas ${\Gamma}_1$ encircles ${\Gamma}_2$ and the
point $z$
(see Figure 4 (a)).

The expansion \expanii\ is convergent for $z$ inside the Cassini annulus
(see Figure 5)
\eqn\domainii{
A_{z_1,z_2}\equiv\lbrace z\in\Omega\setminus\Omega_0, \hskip 2mm
r_2<\vert(z-z_1)(z-z_2)\vert<r_1\rbrace }
where
$$
r_1\equiv {\rm Inf}_{w\in\Css\setminus\Omega}\left\lbrace
\vert(w-z_1)(w-z_2)\vert\right\rbrace,\quad
r_2\equiv{\rm Sup}_{w\in\Omega_0}\left\lbrace
\vert(w-z_1)(w-z_2)\vert\right\rbrace.
$$
}

\noindent{\bf Proof.} By Cauchy's theorem,
\eqn\cauchybis{
f(z)={1\over 2\pi i}\int_{{\Gamma}_1}{f(w)dw\over w-z}-
{1\over 2\pi i}\int_{{\Gamma}_2}{f(w)\,dw\over w-z},
}
where ${\Gamma}_1$ and ${\Gamma}_2$ are the contours defined above. We
substitute
\lequal-\ui\ into the first integral above and
$$
{1\over w-z}={z_1+z_2-z-w\over (z-z_1)(z-z_2)}{1\over 1-u}, \hskip 2cm
u\equiv {(w-z_1)(w-z_2)\over (z-z_1)(z-z_2)},
$$
into the second one. Now we introduce the expansion \expanu\ of the factor
$(1-u)^{-1}$ in both integrals in \cauchybis. After straightforward
calculations we obtain
\expanii-\remii.

For any $z$ verifying \domainii, we can take simple closed loops
${\Gamma}_1$ and
${\Gamma}_2$ in $\Omega\setminus\Omega_0$ such that
$\vert(z-z_1)(z-z_2)\vert<\vert(w-z_1)(w-z_2)\vert$ $\forall$
$w\in{\Gamma}_1$ and
$\vert(z-z_1)(z-z_2)\vert>\vert(w-z_1)(w-z_2)\vert$ $\forall$ $w\in{\Gamma}_2$
(see Figure 4 (b)).
On these contours $\vert f(w)\vert$ is bounded by some constant $C$: $\vert
f(w)\vert\le C$.
Introducing these bounds in \remii\ we see that
$\lim_{N\to\infty}r_N(z_1,z_2;z)=0$ and
the proof follows.
\hfill $\boxe$

\bigskip
\centerline{\epsfxsize=10cm \epsfbox{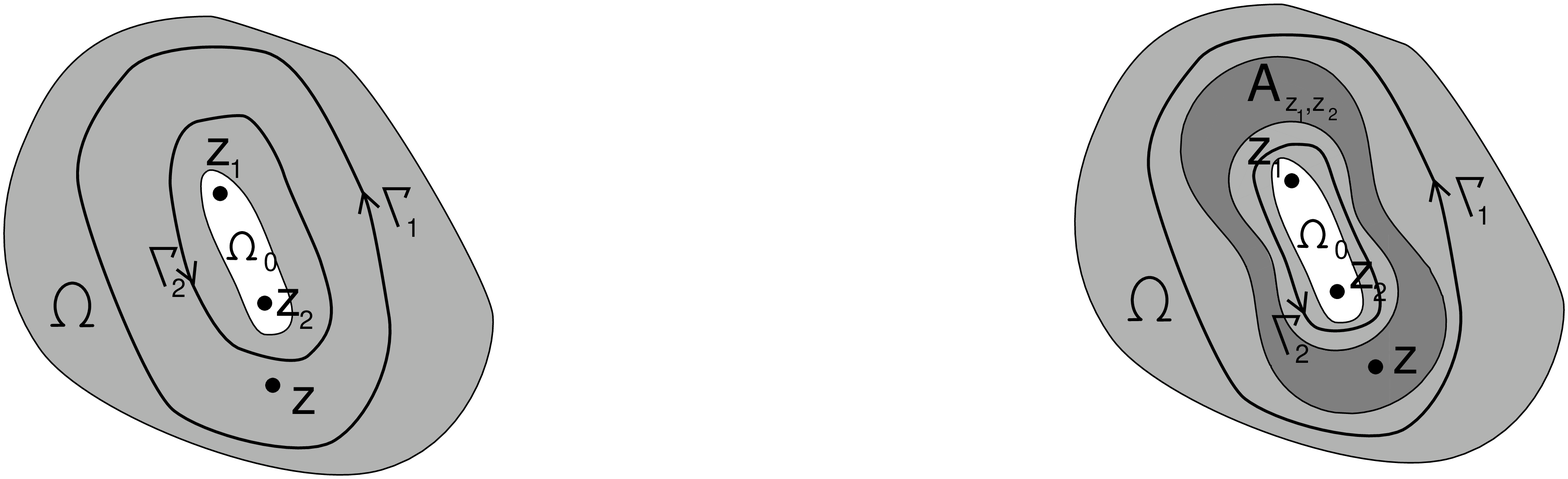}}
\centerline{(a) \hskip 6cm (b)}
\parindent=10pt

\noindent
{\bf Figure 4}. {\cfont (a) Contours $\Gamma_1$ and $\Gamma_2$ in the integrals
\expanii-\remii. (b) For $z\in A_{z_1,z_2}$, we can take a contour
$\Gamma_2$ in
$\Omega$ situated between $\Omega_0$ and $A_{z_1,z_2}$ and a contour
$\Gamma_1$ in $\Omega$ which contains $A_{z_1,z_2}$ inside. Therefore,
$\vert(z-z_1)(z-z_2)\vert<\vert(w-z_1)(w-z_2)\vert$ $\forall$ $w\in\Gamma_1$
and
$\vert(w-z_1)(w-z_2)\vert<\vert(z-z_1)(z-z_2)\vert$ $\forall$
$w\in\Gamma_2$.}
\vskip 2mm

\rm
\parindent=15pt

\bigskip
\centerline{\epsfxsize=10cm \epsfbox{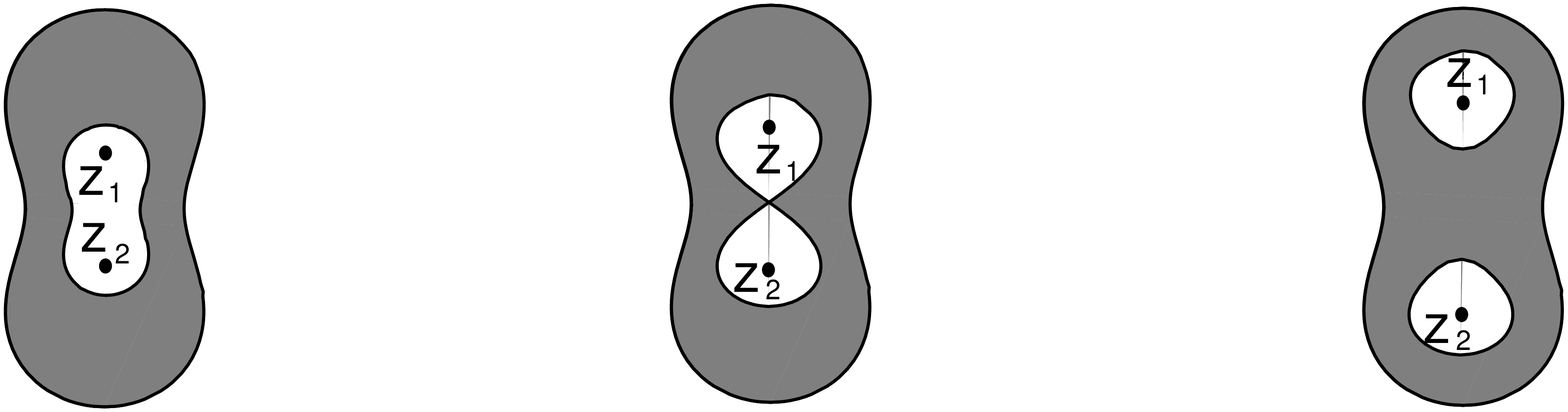}}
\centerline{(a) $4r_1>4r_2>\vert z_1-z_2\vert^2$ \hskip 6mm
(b) $4r_1>\vert z_1-z_2\vert^2=4r_2$
\hskip 6mm (c) $4r_1>\vert z_1-z_2\vert^2>4r_2$}
\bigskip
\centerline{\epsfxsize=7cm \epsfbox{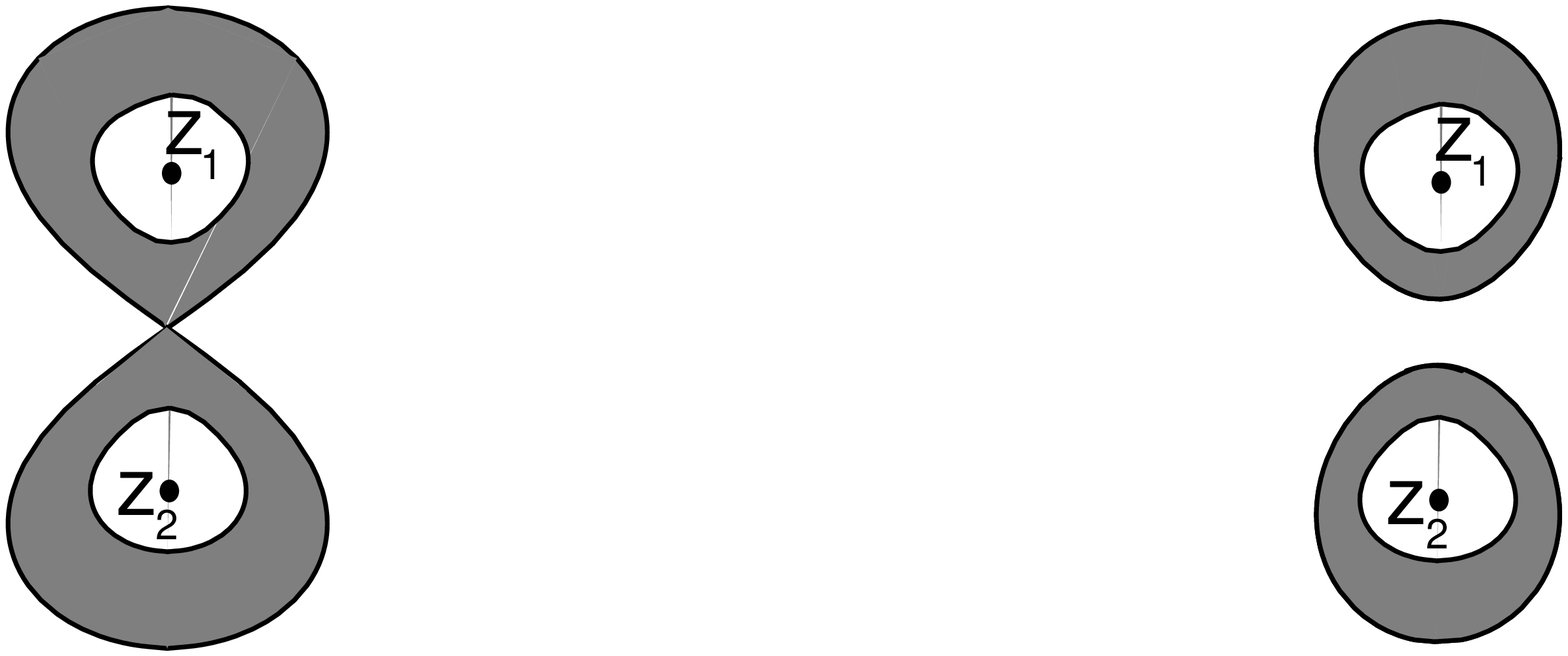}}
\centerline{
(d) $4r_1=\vert z_1-z_2\vert^2>4r_2$
\hskip 2cm (e) $\vert z_1-z_2\vert^2>4r_1>4r_2$}
\parindent=10pt

\noindent
{\bf Figure 5}. {\cfont Shape of the Cassini annulus depending on the relative
size of the parameters $r_1$, $r_2$ and
the focal distance $\vert z_1-z_2\vert$.}
\vskip 2mm

\rm
\parindent=15pt

If the only singularities of $f(z)$ inside $\Omega_0$ are just poles at
$z_1$ and/or
$z_2$, then
an alternative formula to \coefii\ and \coefiii\ to compute the
coefficients of the above two-point
Laurent expansion is given in the following proposition.

\noindent
{\bf Proposition 3.} {\it Suppose that $g_1(z)\equiv (z-z_1)^{m_1}f(z)$
and $g_2(z)\equiv (z-z_2)^{m_2}f(z)$ are
analytic functions in $\Omega$ for certain $m_1$, $m_2\in\Ns$. Then,
for $n=0,1,2,...$,
coefficients $b_n(z_1,z_2)$ and $c_n(z_1,z_2)$ in the expansion \expanii\ are
also given by the formulas:
\eqn\coefiibis{
\eqalign{
b_n(z_1,z_2)= & \sum_{k=0}^{n+m_1-1}\left(\matrix{n+m_1-1 \cr
k}\right){(-1)^{k+1}(n+1)_kg_1^{(n+m_1-k-1)}(z_1)\over
(n+m_1-1)!(z_1-z_2)^{n+k+2}}+ \cr &
\sum_{k=0}^{n+m_2}\left(\matrix{n+m_2 \cr
k}\right){(-1)^{k}(n)_kg_2^{(n+m_2-k)}(z_2)\over
(n+m_2)!(z_2-z_1)^{n+k+1}}, \cr}
}
where $(n)_k$ denotes the Pochhammer symbol and
\eqn\coefiibisbis{
\eqalign{
c_n(z_1,z_2)= & -\sum_{k=0}^{m_1-n-2}k!\left(\matrix{m_1-n-2 \cr
k}\right)\left(\matrix{n \cr k}\right){(z_1-z_2)^{n-k-1}
g_1^{(m_1-n-k-2)}(z_1)\over (m_1-n-2)!}+ \cr &
\sum_{k=0}^{m_2-n-1}k!\left(\matrix{m_2-n-1 \cr
k}\right)\left(\matrix{n+1 \cr k}\right){(z_2-z_1)^{n-k}
g_2^{(m_2-n-k-1)}(z_2)\over (m_2-n-1)!}. \cr}
}
In these formulas, empty sums must be understood as zero. Coefficients
$b_n(z_2,z_1)$ and
$c_n(z_2,z_1)$ are given, respectively, by \coefiibis\ and \coefiibisbis\
interchanging $z_1$, $g_1$ and $m_1$ by $z_2$, $g_2$ and $m_2$ respectively.
}

\noindent
{\bf Proof.} We deform both, the contour ${\Gamma}_1$ in equation \coefii\ and
${\Gamma}_2$ in equation \coefiii, to any
contour  of the form ${\cal C}_1\cup{\cal C}_2$ contained in
$\Omega$, where
${\cal C}_1$ (${\cal C}_2$) is a simple closed loop which
encircles the point $z_1$ ($z_2$) in the
counterclockwise
direction and does not contain the point $z_2$ ($z_1$) inside (see Figure 3
(b)). Then,
$$
\eqalign{
b_n(z_1,z_2)= & {1\over 2\pi i(z_2-z_1)}\left\lbrace\int_{{\cal
C}_1}{g_1(w)\over
(w-z_2)^{n+1}}{dw\over (w-z_1)^{n+m_1}}+ \right. \cr  & \left.
\int_{{\cal C}_2}{g_2(w)\over
(w-z_1)^{n}}{dw\over (w-z_2)^{n+m_2+1}}\right\rbrace = \cr &
{1\over z_2-z_1}\left\lbrace{1\over (n+m_1-1)!}{d^{n+m_1-1}\over dw^{n+m_1-1}}
\left.{g_1(w)\over(w-z_2)^{n+1}}\right\vert_{w=z_1}+
\right. \cr & \left.
{1\over  (n+m_2)!}{d^{n+m_2}\over dw^{n+m_2}}
\left.{g_2(w)\over(w-z_1)^{n}}\right\vert_{w=z_2}\right\rbrace \cr}
$$
and
$$
\eqalign{
& c_n(z_1,z_2)= {1\over 2\pi i(z_2-z_1)}\left\lbrace\int_{{\cal C}_1}
{(w-z_2)^{n}g_1(w)\over (w-z_1)^{m_1-n-1}}dw+
\int_{{\cal C}_2}
{(w-z_1)^{n+1}g_2(w)\over (w-z_2)^{m_2-n}}dw\right\rbrace= \cr &
{1\over z_2-z_1}\left\lbrace\left.{d^{m_1-n-2}\over
dw^{m_1-n-2}}
\left\lbrack{(w-z_2)^{n}g_1(w)\over (m_1-n-2)!}
\right\rbrack\right\vert_{w=z_1}+
\left.{d^{m_2-n-1}\over dw^{m_2-n-1}}
\left\lbrack{(w-z_1)^{n+1}g_2(w)\over  (m_2-n-1)!}
\right\rbrack\right\vert_{w=z_2}\right\rbrace
\cr}
$$
From here, equations  \coefiibis\ and \coefiibisbis\ follow
after straightforward computations.
\hfill$\boxe$

{\bf Remark 2.} {
Let $z$ be a real or complex variable and $z_1, z_2$
($z_1\ne z_2$) two real or
complex numbers. Suppose that $g_1(z)\equiv (z-z_1)^{m_1}f(z)$
is $n-$times differentiable at $z_1$ and  $g_2(z)\equiv (z-z_2)^{m_2}f(z)$
is $n-$times differentiable at $z_2$. Define
$$
g(z)\equiv f(z)-\sum_{n=0}^{M-1}\left\lbrack
c_n(z_1,z_2)(z-z_1)+c_n(z_2,z_1)(z-z_2)\right\rbrack
(z-z_1)^{-n-1}(z-z_2)^{-n-1},
$$
where $M\equiv$Max$\lbrace m_1,m_2\rbrace$. Then, the thesis of Proposition 2
holds for $f(z)$ replaced by $g(z)$. Moreover, if
$(z-z_1)^{m_1}(z-z_2)^{m_2}f(z)$ is an
analytic function in $\Omega$, then the thesis of Theorem 1
applies to $g(z)$.}

\sect{Two-point Taylor-Laurent expansions}

\noindent
{\bf Theorem 3.} {\it Let $\Omega_0$ and $\Omega$ be closed and open
sets, respectively, of the complex
plane, and $\Omega_0\subset\Omega\subset\Cs$. Let $f(z)$
be an analytic function on $\Omega\setminus\Omega_0$,
$z_1\in\Omega_0$ and $z_2\in\Omega\setminus\Omega_0$. Then, for
$z\in\Omega\setminus\Omega_0$, $f(z)$ admits the Taylor-Laurent expansion
\eqn\expaniii{
\eqalign{
f(z)= & \sum_{n=0}^{N-1}\left\lbrack
d_n(z_1,z_2)(z-z_1)+d_n(z_2,z_1)(z-z_2)\right\rbrack
(z-z_1)^n(z-z_2)^n+ \cr &
\sum_{n=0}^{N-1} e_n(z_1,z_2)
(z-z_2)^{n}(z-z_1)^{-n-1}+r_N(z_1,z_2;z), \cr}
}
where the coefficients $d_n(z_1,z_2)$, $d_n(z_2,z_1)$ and $e_n(z_1,z_2)$
of the expansion are given by the  Cauchy integrals
\eqn\coefiv{
d_n(z_1,z_2)\equiv{1\over 2\pi i(z_2-z_1)}\int_{{\Gamma}_1}
{f(w)\,dw\over (w-z_1)^n(w-z_2)^{n+1}}
}
and
\eqn\coefv{
e_n(z_1,z_2)\equiv{z_1-z_2\over 2\pi i}\int_{{\Gamma}_2}
{(w-z_1)^{n}\over(w-z_2)^{n+1}}f(w)\,dw.
}
The remainder term $r_N(z_1,z_2;z)$ is given by the Cauchy integrals
\eqn\remiii{
\eqalign{
r_N(z_1,z_2;z)\equiv & {1\over 2\pi i}\int_{{\Gamma}_1}{f(w)\,dw\over
(w-z_1)^N(w-z_2)^N
(w-z)}(z-z_1)^N(z-z_2)^N- \cr &
{1\over 2\pi i}\int_{{\Gamma}_2}{(w-z_1)^Nf(w)\,dw\over (w-z_2)^N
(w-z)}{(z-z_2)^N\over(z-z_1)^N}.\cr}
}
In these integrals, the contours of integration ${\Gamma}_1$ and
${\Gamma}_2$ are
simple closed loops
contained in $\Omega\setminus\Omega_0$ which encircle $\Omega_0$ in the
counterclockwise
direction.  Moreover, ${\Gamma}_2$ does not contain
the points $z$ and $z_2$ inside,  whereas ${\Gamma}_1$ encircles
${\Gamma}_2$ and the points
$z$ and $z_2$ (see Figure 6 (a)).

The expansion \expaniii\ is convergent in the region (Figure 7)
\eqn\domainiii{
D_{z_1,z_2}\equiv\lbrace z\in\Omega\setminus\Omega_0, \hskip 2mm
\vert(z-z_1)(z-z_2)\vert<r_1
\hskip 2mm {\rm and}\hskip 2mm
\vert z-z_2\vert<r_2\vert z-z_1\vert\rbrace }
where $r_1\equiv$ Inf$_{w\in\Css\setminus\Omega}\left\lbrace
\vert(w-z_1)(w-z_2)\vert\right\rbrace$
and $r_2\equiv$ Inf$_{w\in\Omega_0}\left\lbrace \vert(w-z_2)(w-z_1)^{-1}
\vert\right\rbrace$.
}

\noindent
{\bf Proof.} By Cauchy's theorem,
\eqn\cauchybisbis{
f(z)={1\over 2\pi i}\int_{{\Gamma}_1}{f(w)\,dw\over w-z}-
{1\over 2\pi i}\int_{{\Gamma}_2}{f(w)\,dw\over w-z},
}
where ${\Gamma}_1$ and ${\Gamma}_2$ are the contours defined above. We
substitute \lequal-\ui\ into the first integral above and
\eqn\equalbisbis{
{1\over w-z}={z_2-z_1\over (z-z_1)(w-z_2)}{1\over 1-u}, \hskip 2cm
u\equiv {(w-z_1)(z-z_2)\over (z-z_1)(w-z_2)}
}
 into the second
one. Now we introduce the expansion \expanu\ of the factor $(1-u)^{-1}$ in
both integrals in \cauchybisbis. After straightforward calculations we obtain
\expaniii-\remiii.

For any $z$ verifying \domainiii, we can take simple closed loops
${\Gamma}_1$ and ${\Gamma}_2$ in $\Omega\setminus\Omega_0$ such that
$\vert(z-z_1)(z-z_2)\vert<\vert(w-z_1)(w-z_2)\vert$ $\forall$
$w\in{\Gamma}_1$ and
$\vert(z-z_1)(w-z_2)\vert>\vert(w-z_1)(z-z_2)\vert$ $\forall$ $w\in{\Gamma}_2$
(see Figure 6 (b)).
On these contours $\vert f(w)\vert$ is bounded by some constant $C$: $\vert
f(w)\vert\le C$.
Introducing these bounds in \remiii\ we see that
$\lim_{N\to\infty}r_N(z_1,z_2;z)=0$ and
the proof follows.
\hfill $\boxe$

\bigskip
\centerline{\epsfxsize=10cm \epsfbox{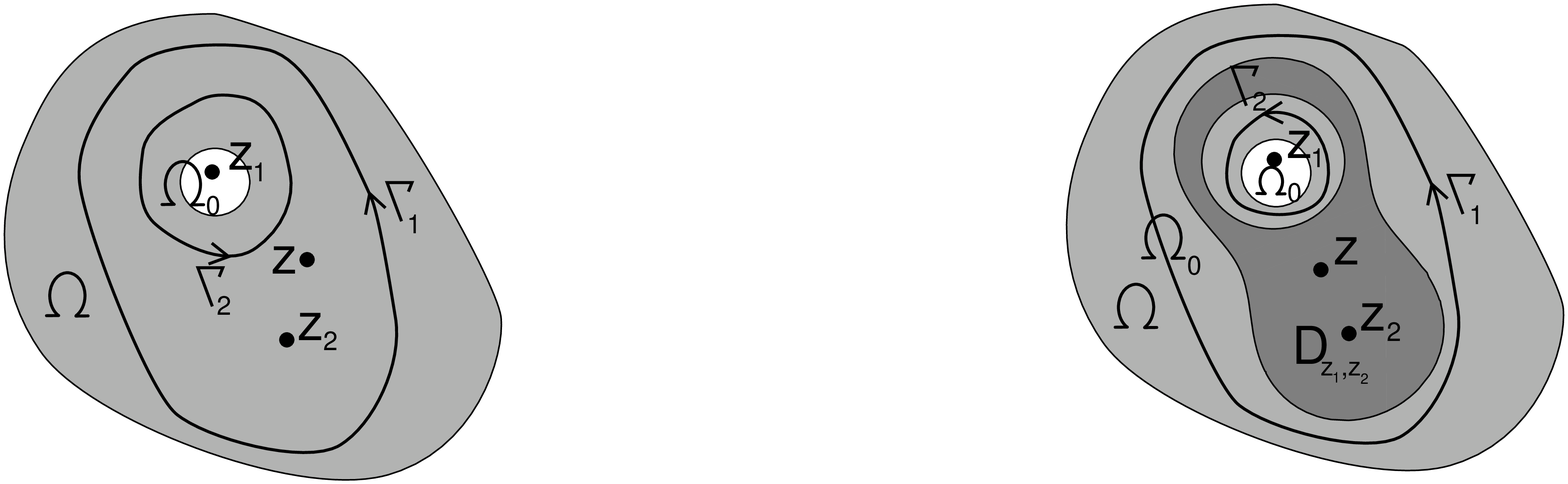}}
\centerline{(a) \hskip 6cm (b)}
\parindent=10pt

\noindent
{\bf Figure 6}. {\cfont (a) Contours $\Gamma_1$ and $\Gamma_2$ in the integrals
\expaniii-\remiii. (b) For $z\in D_{z_1,z_2}$, we can take a contour $\Gamma_2$
situated between $\Omega_0$ and $D_{z_1,z_2}$ and a contour $\Gamma_1$
in $\Omega$ which
contains $D_{z_1,z_2}$ inside. Therefore,
$\vert(z-z_1)(z-z_2)\vert<\vert(w-z_1)(w-z_2)\vert$ $\forall$ $w\in\Gamma_1$
and
$\vert(w-z_1)(z-z_2)\vert<\vert(z-z_1)(w-z_2)\vert$ $\forall$ $w\in\Gamma_2$.}
\vskip 2mm

\rm
\parindent=15pt

\bigskip
\centerline{\epsfxsize=10cm \epsfbox{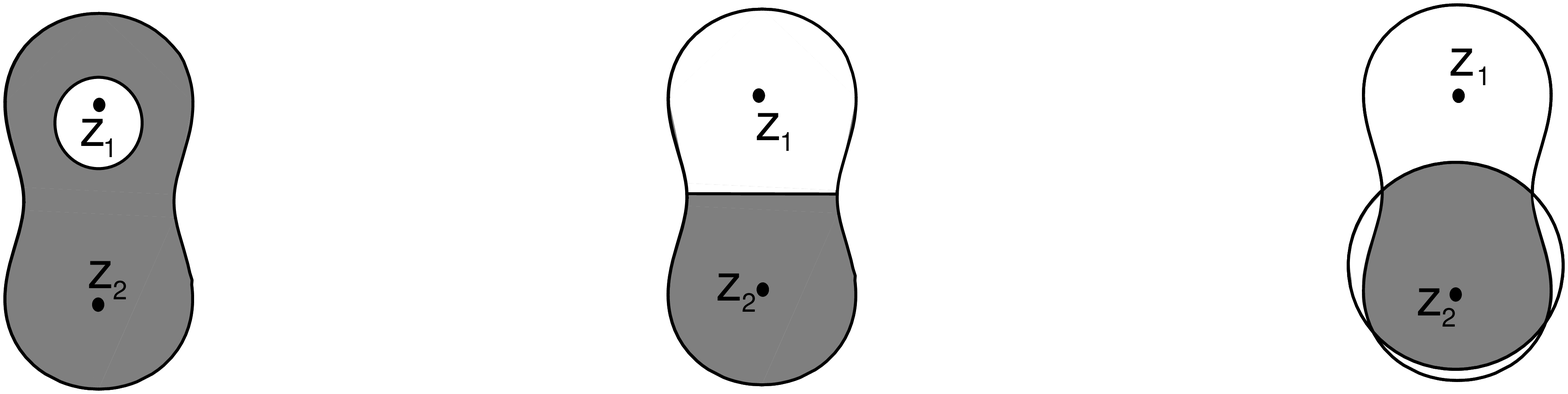}}
\centerline{(a) $4r_1>\vert z_1-z_2\vert^2$, $r_2>1$ \hskip 3mm
(b) $4r_1>\vert z_1-z_2\vert^2$, $r_2=1$
\hskip 3mm (c) $4r_1>\vert z_1-z_2\vert^2$, $r_2<1$}
\bigskip
\centerline{\epsfxsize=10cm \epsfbox{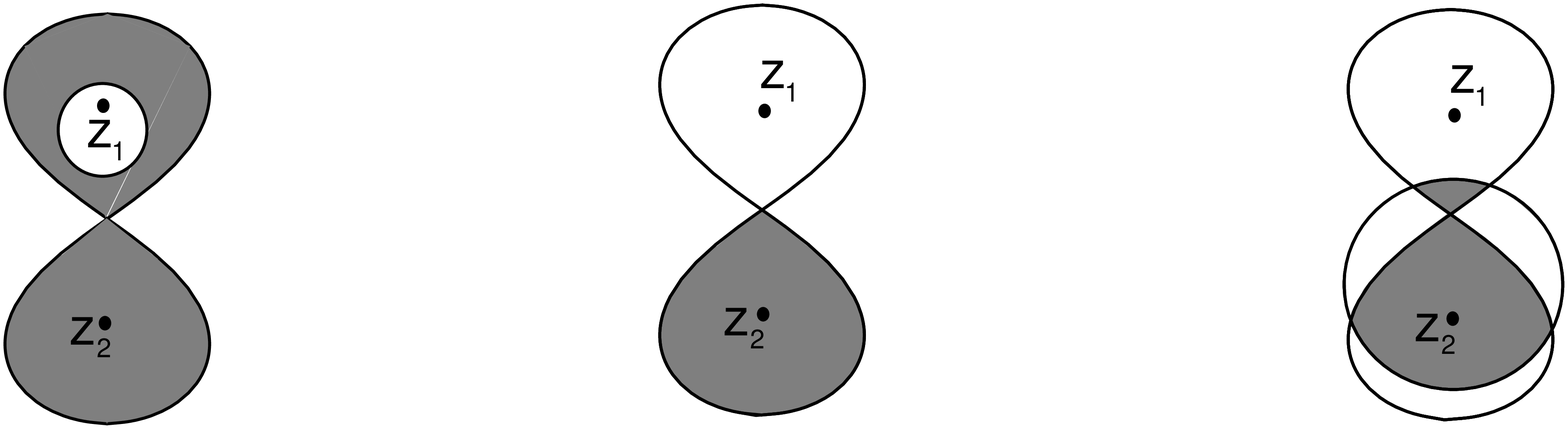}}
\centerline{(d) $4r_1=\vert z_1-z_2\vert^2$, $r_2>1$ \hskip 3mm
(e) $4r_1=\vert z_1-z_2\vert^2$, $r_2=1$
\hskip 3mm (f) $4r_1=\vert z_1-z_2\vert^2$, $r_2<1$}
\bigskip
\centerline{\epsfxsize=11cm \epsfbox{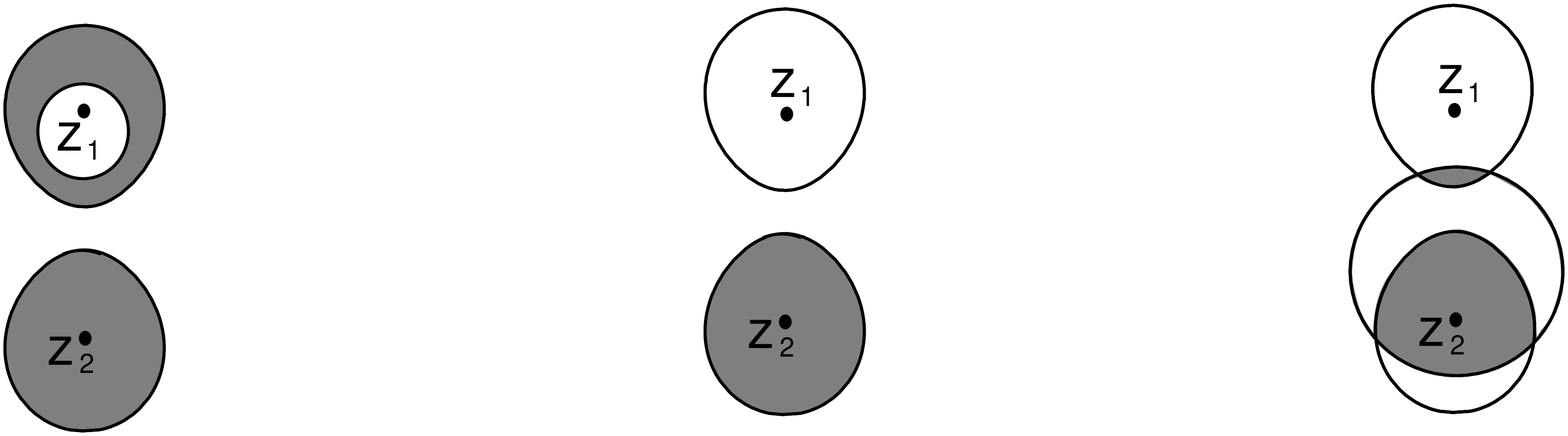}}
\centerline{(g) $4r_1<\vert z_1-z_2\vert^2$, $r_2>1$ \hskip 3mm
(h) $4r_1<\vert z_1-z_2\vert^2$, $r_2=1$
\hskip 3mm (i) $4r_1<\vert z_1-z_2\vert^2$, $r_2<1$}
\parindent=10pt

\noindent
{\bf Figure 7}. {\cfont The region $D_{z_1,z_2}$ defined in Theorem 3 is
given by
$D_{z_1,z_2}=D_1\bigcap D_2$, where $D_1$
is the Cassini oval of focus $z_1$ and $z_2$ and parameter
$r_1$. On the other hand,
for $r_2<1$ ($r_2>1$), $D_2$ is the interior (exterior)
of the circle of center
$z_1+(1-r_2^2)^{-1}(z_2-z_1)=z_2+r_2^2(r_2^2-1)^{-1}(z_1-z_2)$
and radius $\vert z_1-z_2\vert r_2/\vert r_2^2-1\vert$.
For $r_2=1$, $D_2$ is just the half plane
$\vert z-z_2\vert<\vert z-z_1\vert$. The
shape of the Cassini annulus depends on the relative
size of the parameters $\sqrt{r_1}$, $\sqrt{r_2}$ and
the focal distance $\vert z_1-z_2\vert$.}
\vskip 2mm

\rm
\parindent=15pt

If the only singularities of $f(z)$ inside $\Omega_0$ are just poles at
$z_1$, then
an alternative formula to \coefiv-\coefv\ to compute the coefficients of
the above two-point
Taylor-Laurent expansion is given in the following proposition.

\noindent
{\bf Proposition 4.} {\it Suppose that $g(z)\equiv (z-z_1)^{m}f(z)$ is an
analytic function in $\Omega$ for certain $m\in\Ns$. Then, coefficients
$d_n(z_1,z_2)$ and $d_n(z_2,z_1)$ in the expansion \expaniii\ are also
given by the formulas:
\eqn\do{
d_0(z_1,z_2)=  {f(z_2)\over z_2-z_1}-\sum_{k=0}^{m-1}
{1\over (m-k-1)!}{g^{(m-k-1)}(z_1)\over(z_2-z_1)^{k+2}}, \hskip 1cm
}
$$
d_0(z_2,z_1)=  {1\over m!}{g^{(m)}(z_1)\over z_1-z_2},
$$
and, for $n=1,2,3...$,
\eqn\coefiiibis{
\eqalign{
d_n(z_1,z_2)= & -{(-1)^n\over n!}\left\lbrace \sum_{k=0}^{m+n-1}
{(n+k)!\over k!(m+n-k-1)!}{g^{(m+n-k-1)}(z_1)\over(z_2-z_1)^{n+k+2}}+
\right. \cr &
\left. n\sum_{k=0}^{n}
{(n+k-1)!\over
k!(n-k)!}{f^{(n-k)}(z_2)\over(z_1-z_2)^{n+k+1}}\right\rbrace, \cr}
}
\eqn\coefiiirebis{
\eqalign{
d_n(z_2,z_1)= & -{(-1)^n\over n!}\left\lbrace n\sum_{k=0}^{m+n}
{(n+k-1)!\over k!(m+n-k)!}{g^{(m+n-k)}(z_1)\over(z_2-z_1)^{n+k+1}}+ \right.
\cr &
\left.\sum_{k=0}^{n-1}
{(n+k)!\over
k!(n-k-1)!}{f^{(n-k-1)}(z_2)\over(z_1-z_2)^{n+k+2}}\right\rbrace.\cr}
}
For $n=0,1,2,...$, coefficients $e_n(z_1,z_2)$ are given by
\eqn\coefiiibisbis{
e_n(z_1,z_2)=  {(-1)^n\over n!}\sum_{k=0}^{m-n-1}
{(n+k)!\over k!(m-n-k-1)!}{g^{(m-n-k-1)}(z_1)\over(z_2-z_1)^{n+k}}.
}
}

\noindent
{\bf Proof.} We deform both, the contour ${\Gamma}_1$ in equation \coefiv\ and
the  contour ${\Gamma}_2$ in equation \coefv\  to any
contour  of the form ${{\cal C}_1}\cup{{\cal C}_2}$ contained in
$\Omega$, where ${{\cal C}_1}$ (${{\cal C}_2}$) is a simple closed loop
which encircles the point $z_1$ ($z_2$) in the counterclockwise
direction and does not contain the point $z_2$ ($z_1$) inside (see Figure 3
(c)). Then,
$$
\eqalign{
d_n(z_1,z_2)= & {1\over 2\pi i(z_2-z_1)}\left\lbrace\int_{{\cal C}_1}{g(w)\over
(w-z_2)^{n+1}}{dw\over (w-z_1)^{n+m}}+
\right. \cr & \left. \int_{{\cal C}_2}{f(w)\over
(w-z_1)^{n}}{dw\over (w-z_2)^{n+1}}\right\rbrace= \cr &
{1\over (z_2-z_1)}\left\lbrace{1\over (n+m-1)!}{d^{n+m-1}\over dw^{n+m-1}}
\left.{g_1(w)\over(w-z_2)^{n+1}}\right\vert_{w=z_1}+
\right. \cr & \left.
{1\over  n!}{d^{n}\over dw^{n}}
\left.{f(w)\over(w-z_1)^{n}}\right\vert_{w=z_2}\right\rbrace, \cr}
$$
an analog formula for $d_n(z_2,z_1)$, and
$$
\eqalign{
e_n(z_1,z_2)= & {z_1-z_2\over 2\pi i}\int_{{\cal C}_1}
{g(w)\over(w-z_2)^{n+1}}{dw\over (w-z_1)^{m-n}}= \cr &
(z_1-z_2){1\over (m-n-1)!}{d^{m-n-1}\over dw^{m-n-1}}
\left. {g(w)\over (w-z_2)^{n+1}}\right\vert_{w=z_1}. \cr}
$$
From here, equations  \do-\coefiiibisbis\ follow
after straightforward computations.
\hfill$\boxe$

{\bf Remark 3.} {
Let $z$ be a real or complex variable and
$z_1$ and $z_2$ ($z_1\ne z_2$) two real or
complex numbers. Suppose that $(z-z_1)^{m}f(z)$
is $n-$times differentiable at $z_1$ for certain $m\in\Ns$ and $f(z)$
is $n-$times differentiable at $z_2$. Define
$$
g(z)\equiv f(z)-\sum_{n=0}^{m-1}
e_n(z_1,z_2)(z-z_1)^{-n-1}(z-z_2)^{n}.
$$
Then, the thesis of Proposition 2
holds for $g(z)$. If moreover, $(z-z_1)^{m}f(z)$ is an
analytic function in $\Omega$, then the thesis of Theorem 1
applies to $g(z)$.}

\sect{Acknowledgements}

\noindent
J. L. L\'opez wants to thank the C.W.I. of Amsterdam for its scientific and
financial support during the realization of this work.  The financial
support of the saving bank {\it Caja Rural de Navarra} is also acknowledged.

\listrefs

\end